\documentclass{amsart}

\usepackage{amsmath}
\usepackage{amscd}
\usepackage{amssymb}

\newcommand{\cal}{\mathcal}

\newcommand{\bC}{{\mathbb C}}
\newcommand{\bP}{{\mathbb P}}
\newcommand{\bQ}{{\mathbb Q}}
\newcommand{\bR}{{\mathbb R}}
\newcommand{\bZ}{{\mathbb Z}}

\newcommand{\cO}{{\cal O}}

\DeclareMathOperator{\Img}{Im} \DeclareMathOperator{\ch}{ch}
\DeclareMathOperator{\td}{td} \DeclareMathOperator{\Res}{Res}

\DeclareMathOperator{\Ell}{Ell}\DeclareMathOperator{\Ind}{Ind}

\newtheorem{theorem}{Theorem}[section]
\newtheorem{theorem/definition}{Theorem/Definition}[section]

\newtheorem{proposition}{Proposition}[section]
\newtheorem{lemma}{Lemma}[section]

\newtheorem{corollary}{Corollary}[section]

\theoremstyle{remark}
\newtheorem{remark}{Remark}[section]

\theoremstyle{definition}
 
\newtheorem{definition}{Definition}[section]

\begin{document}
\title{Elliptic Genera of Complete Intersections}
\author{Xiaoguang Ma and Jian Zhou}
\address{Department of Mathematical Sciences\\
Tsinghua University\\
Beijing 100084 \\
People's Republic of China} \email{jzhou@math.tsinghua.edu.cn}

\begin{abstract}
We propose a new definition of the elliptic genera for complete intersections,
not necessarily nonsingular,
in projective spaces.
We also prove they coincide with the expressions obtained from Landau-Ginzburg model
by an elementary argument.
\end{abstract}
\date{}
\maketitle

\section{Introduction}

In this work we give a definition of the two-variable elliptic genera
of complete intersections.
In the case of hypersurfaces,
we prove the elliptic genera coincide with the elliptic genera of orbifold
Landau-Ginzburg  model proposed in \cite{Egu-Jin}.
We also extend our results to the case of complete intersections
which has not been treated in the literature.

The problem of defining elliptic genera for singular varieties
has been addressed in  \cite{Bor-Lib},
where the definitions for orbifolds,
for $\bQ$-Gorenstein varieties with log terminal singularities,
and for pairs consisting of projective varieties
and $\bQ$-Cartier divisor on them are given.
The authors of that work used the resolutions of the singularities and
proved their definition is independent of the resolutions.
In this paper,
we take a completely different approach.
More precisely, we will use the resolution of the structure sheaf,
instead of the variety itself.

We will also use an elementary method to establish the relationship between elliptic genera of
complete intersections and LG orbifolds.
We first express the elliptic genera as residues by Hirzebruch-Riemann-Roch theorem,
then apply the residue theorem in one complex variable.
A completely different method has been used in \cite{Gor-Mal}
to treat the case of Calabi-Yau hypersurfaces.

The rest of this paper is arranged as follows.
We recall in Section \ref{sec:Pre} some preliminaries
such as definitions and the Hirzebruch-Riemann-Roch theorem.
We also prove two simple lemmas on residues of analytic functions.
In Section \ref{sec:Smooth} we first
study the Euler numbers, elliptic genera, and $\chi_y$-genera
of smooth hypersurfaces and complete intersections,
then make the definition in the singular case.
In Section \ref{sec:Projective},
we express Euler numbers, elliptic genera, and $\chi_y$-genera
of hypersurfaces and complete intersections in projective spaces
first as residues,
then compute their generating functions.
The results for Euler numbers and the $\chi_y$-genera are classically
known \cite{Hir}.
In Section \ref{sec:Wit} we make generalizations to Witten genera.
Finally in Section \ref{sec:LG} we establish the relationship with Landau-Ginzburg model.

{\bf Acknowledgements}.
{\em The research of the second author is supported by grants from NSFC and Tsinghua University.
He thanks Kefeng Liu for helpful discussions and suggestions.}

\section{Preliminaries}
\label{sec:Pre}

\subsection{Some multiplicative operations on graded vector spaces}
Let $E$ be a finite dimensional vector space over $\bC$.
Consider the formal power series
\begin{align*}
\Lambda_t(E) & = \sum_{n \geq 0} t^n \Lambda^nE, &
S_t(E) & = \sum_{n \geq 0} t^n S^n E.
\end{align*}
The following properties are well-known:
\begin{align*}
& \Lambda_t(E \oplus \widetilde{E}) = \Lambda_t(E) \Lambda_t(\widetilde{E}), \\
& S_t(E \oplus \widetilde{E}) = S_t(E) S_t(\widetilde{E}), \\
& \Lambda_t(E)S_{-t}(E) = 1.
\end{align*}
We understand them as identities of formal power series with coefficients
in the Grothedieck ring $K(*)$ of the isomorphism classes of
finite dimensional vector spaces over $\bC$.
An element of this ring can be written as $E_0 - E_1$ for two finite dimensional
vector spaces $E_0$ and $E_1$.
As usual,
we define
\begin{align*}
\Lambda_t(E_0 - E_1) & = \Lambda_t(E_0)S_{-t}(E_1)
= \Lambda_t(E_0)\Lambda_t(E_1)^{-1}, \\
S_t(E_0 - E_1) & = S_t(E_0)\Lambda_{-t}(E_1) = S_t(E_0)S_t(E_1)^{-1}.
\end{align*}
We interpret  $E_0 - E_1$ as the graded space $E_0 \oplus E_1$.
In general,
for a graded vector space $E = E_0 \oplus \dots \oplus E_n$,
we define
\begin{align*}
\Lambda_t(E) & = \Lambda_t(E_0) \Lambda_t(E_1)^{-1} \cdots \Lambda_t(E_n)^{(-1)^n}, \\
S_t(E) & = S_t(E_0) S_t(E_1)^{-1} \cdots S_t(E_n)^{(-1)^n}.
\end{align*}
We regard them as the {\em graded} exterior and symmetric powers of $E$.

\subsection{The $N=2$ superconformal power}

Let $T$ be a finite dimensional $\bZ$-graded complex vector space.
Consider the formal power series of vector spaces
$$\Ell(T; q, y) = y^{-\dim T/2}\otimes_{n \geq 1}
\left(\Lambda_{-yq^{n-1}}(T^*) \otimes \Lambda_{-y^{-1}q^{n}}(T)
\otimes S_{q^n}(T^*) \otimes S_{q^n}(T)\right).$$
We call this the {\em $N=2$ superconformal power of $T$} for the following reason.
Suppose
$$\Ell(T; q, y) = \sum_{m,n} y^mq^nE_{mn},$$
then the graded space $\Ell(T) = \sum_{m,n}E_{mn}$
has a structure of an $N=2$ superconformal vertex algebra (SCVA).
It is straightforward to see that the operation $\Ell$ is multiplicative
in the following sense:
$$\Ell(T_1 \oplus T_2; q, y) = \Ell(T_1;q,y)\otimes \Ell(T_2; q,y).$$

\subsection{Riemann-Roch number}

Suppose $M$ is a compact complex manifold,
and $\pi:E \to M$ is a holomorphic vector bundle on $M$.
The Riemann-Roch number of $E$ is defined by
$$\chi(M, E) = \sum_{k \geq 0} (-1)^k \dim H^k(M, \cO_M(E))$$
where $\cO_M(E)$ is the sheaf of holomorphic section to $E$,
and $H^k(M, \cO_M(E))$ is the $k$-th cohomology of $\cO_M(E)$.
It has the following property:
If
$$0 \to E_1 \to E_2 \to E_3 \to 0$$
is an exact sequence of holomorphic vector bundles on $M$,
then we have
$$\chi(M, E_1) - \chi(M, E_2) + \chi(M, E_3) = 0.$$
Hence if one denotes by $K(M)$ the Grothendieck ring of holomorphic
vector bundles over $M$,
then one gets an additive homomorphism
$$\chi(M, \cdot): K(M) \to \bZ.$$
It is straightforward to extend $\chi(M, \cdot)$ to formal power series
with coefficients in $K(M)$.

The Riemann-Roch number of $E$ can be computed by the Hirzebruch-Riemann-Roch (HRR)
Theorem as follows.
Let $x_1, \dots, x_d$ and $z_1, \dots, z_r$ be the formal Chern roots of $TX$
and $E$ respectively,
then
$$\chi(X, E) = \int_X \ch(E)\td(TX)
= \int_X \sum_j e^{z_j}\prod_i \frac{x_i}{1 - e^{-x_i}}.$$

\subsection{Elliptic genus of a complex manifold}

Suppose $X$ is a compact complex manifold of dimension $d$.
Let $TX$ denote the holomorphic tangent bundle of $X$,
and consider the holomorphic vector bundle $\Ell(TX; q, y)$.

\begin{definition}
The formal power series
\begin{eqnarray*}
\chi(X; q, y) = \chi(X, \Ell(TX; q, y))
\end{eqnarray*}
is called the {\em elliptic genus of $X$}.
\end{definition}

BY HRR, one can easily show that:
\begin{eqnarray} \label{eqn:chiX}
&& \chi(X, q, y) = y^{-d/2} \int_X \prod_{j=1}^d x_j
\prod_{n \geq 1}\frac{(1-yq^{n-1}e^{-x_j})(1- y^{-1}q^ne^{x_j})}
{(1-q^{n-1}e^{-x_j})(1-q^ne^{x_j})}.
\end{eqnarray}

\subsection{Some preliminary results on residues}
For later applications,
we prove here some results on residues of a sequence of functions.

\begin{lemma} \label{lm:Residue}
Let $f(z)$ be analytic near $z=0$,
$f(0) = 0$,
and $f'(0) \neq 0$.
For $t$ close to $0$,
let $z = g(t)$ satisfies $f(g(t)) = t$.
Then for any constant $m \neq 0$ we have
\begin{eqnarray*}
&& \sum_{n =1}^{\infty} t^n \Res_z \frac{f(mz)}{f(z)^{n+1}}
= \frac{f(mg(t))}{f'(g(t))}.
\end{eqnarray*}
More generally,
for nonzero constants $m_1, \dots, m_r$,
we have
\begin{eqnarray*}
&& \sum_{n = r}^{\infty} t^n \Res_z \frac{f(m_1z)\cdots
f(m_rz)}{f(z)^{n+1}} = \frac{f(m_1g(t))\cdots
f(m_rg(t))}{f'(g(t))}.
\end{eqnarray*}
\end{lemma}

\begin{proof}
Fix a small enough $\epsilon > 0$.
We have
\begin{eqnarray*}
\sum_{n =1}^{\infty}t^n \Res_z \frac{f(mz)}{f(z)^{n+1}} &  = &
\frac{1}{2\pi \sqrt{-1}} \sum_{n=1}^{\infty} \int_{|z| = \epsilon}
\frac{f(mz)}{f(z)^{n+1}}t^n dz \\
& = & \frac{1}{2\pi \sqrt{-1}} \int_{|z| = \epsilon}
\sum_{n=1}^{\infty}
\frac{f(mz)}{f(z)^{n+1}}t^n dz \\
& = & \frac{1}{2\pi \sqrt{-1}} \int_{|z| = \epsilon}
\frac{tf(mz)}{f(z)^2}\frac{1}{1 - \frac{t}{f(z)}} dz \\
& = & \frac{1}{2\pi \sqrt{-1}} \int_{|z| = \epsilon}
\frac{tf(mz)}{f(z)}\frac{1}{f(z) - t} dz.
\end{eqnarray*}
Now the integrand in the last equality has only a simple pole at
$z = g(t)$ on the disk $|z| \leq \epsilon$,
hence by residue theorem,
\begin{eqnarray*}
\sum_{n =1}^{\infty}t^n \Res_z \frac{f(mz)}{f(z)^{n+1}}
&  = & \left. \frac{tf(mz)}{f(z)}\frac{1}{\frac{d}{dz}(f(z) - t)}\right|_{z = g(t)} \\
& = & \frac{tf(mg(t))}{f(g(t))}\frac{1}{f'(g(t))} \\
& = & \frac{f(mg(t))}{f'(g(t))}.
\end{eqnarray*}
The general case is proved in the same fashion.
\end{proof}

Sometimes it is not easy to find $g(t)$ explicitly.
We offer the following calculations for small $n$ directly from $f$.

\begin{lemma} \label{lm:Residue2}
Suppose we are in the situation of Lemma \ref{lm:Residue}.
Then we have
\begin{eqnarray*}
\Res_z \frac{f(mz)}{f(z)^{n+1}}
= \begin{cases}
\frac{m}{f'(0)}, & n = 1, \\
\frac{m(m-3)}{2} \frac{f''(0)}{f'(0)^3}, & n =2, \\
\frac{m^3-4m}{6}\cdot \frac{f'''(0)}{f'(0)^4} - \frac{2m^2 -
5m}{2} \cdot\frac{f''(0)^2}{f'(0)^5}, & n = 3.
\end{cases}
\end{eqnarray*}
\end{lemma}

\begin{proof}
We have
\begin{eqnarray*}
f(z) & = & f'(0)z + \frac{f''(0)}{2}z^2 + \frac{f'''(0)}{6}z^3 + \dots
\end{eqnarray*}
Hence
\begin{eqnarray*}
\frac{f(mz)}{f(z)^{n+1}}
& = & \frac{f'(0)mz + \frac{f''(0)}{2}m^2z^2 + \frac{f'''(0)}{6}m^3z^3 + \dots}
{(f'(0)z + \frac{f''(0)}{2}z^2 + \frac{f'''(0)}{6}z^3 + \dots)^{n+1}} \\
& = & \frac{m}{(f'(0))^{n}z^n} \frac{1 + \frac{f''(0)}{2f'(0)}mz +
\frac{f'''(0)}{6f'(0)}m^2z^2 + \dots} {(1 + \frac{f''(0)}{2f'(0)}z
+ \frac{f'''(0)}{6f'(0)}z^2 + \dots)^{n+1}}
\end{eqnarray*}
The Lemma is proved by Laurent expansion.
\end{proof}

\subsection{Theta functions}
\label{sec:theta}

In this subsection we collect the definitions and properties of
the theta functions which will be used later.
For $\tau \in \bC$ with $\Img \tau > 0$,
define the theta functions as four infinite products \cite{Cha}:
\begin{eqnarray*}
\theta(\upsilon,\tau)
& = & \frac{\sqrt{-1}q^{\frac{1}{8}}}{e^{\pi\sqrt{-1}\upsilon}}
\prod_{k\geq1} [(1-q^k) (1-q^{k-1}e^{2\pi \sqrt{-1}\upsilon})
(1-q^ke^{-2\pi \sqrt{-1} \upsilon})], \\
\theta_{1}(\upsilon,\tau) & = & q^{\frac{1}{8}}e^{\pi \sqrt{-1} \upsilon}
\prod_{k\geq1} [(1-q^k)(1+q^ke^{2\pi \sqrt{-1} \upsilon})
(1+q^{k-1}e^{-2\pi \sqrt{-1} \upsilon})],\\
\theta_{2}(\upsilon,\tau)& = & \prod_{k\geq1}
[(1-q^k)(1-q^{k-\frac{1}{2}}e^{2\pi \sqrt{-1} \upsilon})
(1-q^{k-\frac{1}{2}}e^{-2\pi \sqrt{-1} \upsilon})], \\
\theta_{3}(\upsilon,\tau)& =
&\prod_{k\geq1} [(1-q^k)(1+q^{k-\frac{1}{2}}e^{2\pi \sqrt{-1} \upsilon})
(1+q^{k-\frac{1}{2}}e^{-2\pi\sqrt{-1} \upsilon})],
\end{eqnarray*}
where $q=e^{2\pi\sqrt{-1}\tau}$.

Recall the following properties of the theta functions.
\begin{align*}
\theta(\upsilon+1, \tau)  & = -\theta(\upsilon,\tau), &
\theta(\upsilon+\tau, \tau) & =-q^{-\frac{1}{2}}e^{-2\pi\sqrt{-1}\upsilon}\theta(\upsilon,\tau), \\
\theta_{1}(\upsilon+1,\tau)  & =-\theta_{1}(\upsilon,\tau), &
\theta_{1}(\upsilon+\tau,\tau)
& =q^{-\frac{1}{2}}e^{-2\pi\sqrt{-1}\upsilon}\theta_{1}(\upsilon,\tau), \\
\theta_{2}(\upsilon+1,\tau) & = \theta_{2}(\upsilon,\tau), &
\theta_{2}(\upsilon+\tau,\tau)
& =-q^{-\frac{1}{2}}e^{-2\pi\sqrt{-1}\upsilon}\theta_{2}(\upsilon,\tau),\\
\theta_{3}(\upsilon+1,\tau) & =\theta_{3}(\upsilon,\tau), &
\theta_{3}(\upsilon+\tau,\tau)  & =
q^{-\frac{1}{2}}e^{-2\pi\sqrt{-1}\upsilon}\theta_{3}(\upsilon,\tau)
\end{align*}
The theta functions have the following simple zeros:
\begin{align*}
\theta(a+b\tau,\tau) & = 0, &
\theta_{1}(a+b\tau+\frac{1}{2},\tau) & =0, \\
\theta_{2}(a+b\tau+\frac{\tau}{2},\tau) &=0, &
\theta_{3}(a+b\tau+\frac{1+\tau}{2},\tau) &=0,
\end{align*}
for $a,b\in \bZ$.
One also has:
\begin{align}
\theta_{1}(\upsilon,\tau)  & = \theta(\frac{1}{2}-\upsilon,\tau), \label{eqn:theta1} \\
\theta_{2}(\upsilon,\tau)
& =  -\sqrt{-1}q^{\frac{1}{8}}e^{-\pi
\sqrt{-1}\upsilon}\theta(\frac{\tau}{2}-\upsilon,\tau),  \label{eqn:theta2} \\
\theta_{3}(\upsilon,\tau) & = q^{\frac{1}{8}}e^{-\pi
\sqrt{-1}\upsilon}\theta(\frac{1}{2}+\frac{\tau}{2}-\upsilon,\tau).  \label{eqn:theta3}
\end{align}
In particular,
\begin{align}
& \theta_3(\upsilon, \tau) = \theta_2(\upsilon + \frac{1}{2}, \tau).
\end{align}
Finally,
\begin{align} \label{eqn:-theta}
\theta(-\upsilon,\tau) & = -\theta(\upsilon,\tau), &
\theta_i(-\upsilon,\tau) & = \theta_i(\upsilon,\tau), \qquad i =1,2,3.
\end{align}

\section{Elliptic Genera of Hypersurfaces and Complete Intersections}

\label{sec:Smooth}

In this section we first study the elliptic genera of smooth
hypersurfaces and complete intersections by adjunction formula,
then we reverse the procedure to define elliptic genera
for general complete intersections.

\subsection{Smooth hypersurfaces and adjunction formula}

Let $X$ be a complex manifold and let $Y \subset X$ be a smooth hypersurface.
Denote by $[Y]$ the line bundle on $X$ associated to the divisor $Y$.
The adjunction formula states
$$N_{Y/X} \cong [Y]|_Y,$$
where $N_{Y/X} = TX|_Y/TY$ is the normal bundle of $Y$ in $X$.
Therefore the exact sequence
$$0 \to TY \to TX|_Y \to N_{Y/X} \to 0$$
can be rewritten as
\begin{eqnarray} \label{eqn:Adjunction}
0 \to TY \to TX|_Y \to [Y]|_Y \to 0.
\end{eqnarray}
In $K$-theory,
$$TY = (TX - [Y])|_Y,$$
we refer to the latter as the {\em virtual tangent bundle of $Y$ in $X$}.

\subsection{Euler number of smooth hypersurfaces}

As a warmup exercise,
let us recall the computation of the Euler number of $Y$ by (\ref{eqn:Adjunction}).
First of all,
$$c(TY) = c(TX)c([Y])^{-1}|_Y.$$
We introduce the following notation:
$$c(E) = \prod_i c(E_i)^{(-1)^i}$$
for a $\bZ$-graded vector bundle $E = \oplus_i E_i$.
Then
$$c(TY) = c(TX - [Y])|_Y.$$
Hence
$$\chi(Y) = \int_Y c(TY)
= \int_X c(TY) c_1(N_{Y/X})
= \int_X c(TX - [Y]) c_1([Y]).$$
Thus we have prove the following

\begin{theorem}
Suppose  $Y$ is a smooth hypersurface of a compact complex manifold $X$.
Then we have
\begin{eqnarray*}
&& \chi(Y)
= \int_X c(TX - [Y]) c_1([Y]).
\end{eqnarray*}
\end{theorem}

In terms of formal Chern roots,
we have
\begin{eqnarray} \label{eqn:chi}
\chi(Y) = \int_X
\frac{c\prod_{j=1}^{\dim X} (1+x_j)}{1+c},
\end{eqnarray}
where $c=c_1([Y])$, and $x_i$ are the formal Chern roots of $TX$.

\subsection{Elliptic genus of a smooth hypersurface}

By the multiplicative property of $\Ell$,
one gets from (\ref{eqn:Adjunction}):
$$\Ell(TY;q,y)\Ell([Y]|_Y;q;y) = \Ell(TX|_Y;q,y) = \Ell(TX;q,y)|_Y,$$
or equivalently,
$$\Ell(TY;q,y) = \Ell(TX - [Y];q,y)|_Y.$$
Hence
$$\chi(Y, \Ell(TY;q,y)) = \chi(Y, \Ell(TX - [Y];q,y)|_Y).$$
On the other hand,
as usual, one has the following exact sequence
$$0 \to \cO_X([-Y]) \to \cO_X \to \cO_Y \to 0.$$
Tensoring with formal series of vector bundles $\Ell(TX;q,
y)\Ell([Y];q;y)^{-1}$,
\begin{eqnarray*}
0 & \to & \cO_X(\Ell(TX-[Y];q,y) \to \cO_X(\Ell(TX-[Y];q,y) \\
& \to & \cO_X(\Ell(TX-[Y];q,y)|_Y \to 0,
\end{eqnarray*}
hence by taking the Riemann-Roch number,
one gets
\begin{eqnarray*}
&&\chi(Y, \Ell(TX-[Y];q,y)|_Y) = \chi(X, \Ell(TX-[Y]; q, y)(\cO_X
- [-Y])).
\end{eqnarray*}
Therefore,
we have proved the following

\begin{theorem} \label{tom:Hypersurface}
Suppose  $Y$ is a smooth hypersurface of a compact complex manifold $X$.
Then we have
\begin{eqnarray*}
&& \chi(Y; q, y) =  \chi(X, \Ell(TX-[Y]; q, y)(\cO_X - [-Y])).
\end{eqnarray*}
\end{theorem}

Denote by $c$ the first Chern class of $[Y]$.
Then in terms of the formal Chern roots,
we have by HRR
\begin{eqnarray*}
\chi(Y; q, y)
& = & y^{-\dim Y/2} \int_X \prod_{n \geq 1}
\frac{(1-q^{n-1}e^{-c})(1-q^ne^{c})}{(1-yq^{n-1}e^{-c})(1-y^{-1}q^ne^c)} \\
&& \cdot \prod_{j=1}^{\dim X} x_j \prod_{n \geq 1}
\frac{(1-yq^{n-1}e^{-x_j})(1- y^{-1}q^ne^{x_j})}
{(1-q^{n-1}e^{-x_j})(1-q^ne^{x_j})}.
\end{eqnarray*}

\subsection{The case of Hirzebruch $\chi_y$-genus}
Now we take $q=0$.
Noticing
$$\Ell(TX; 0, y) = y^{-\dim X/2}\Lambda_{-y}(TX)$$ and
$$\Ell([Y]; 0, y) = y^{-1/2}(1 - y [-Y]),$$
one easily gets the following corollary:

\begin{theorem}
Suppose  $Y$ is a smooth hypersurface of a compact complex manifold $X$.
Then we have
\begin{eqnarray*}
&& \chi_{-y}(Y)
= \chi\left(X, \Lambda_{-y}(TX)(1 - y[-Y])^{-1}(1 - [-Y])\right).
\end{eqnarray*}
\end{theorem}

In terms of formal Chern roots,
we have
\begin{eqnarray} \label{eqn:chiy}
\chi_{-y}(Y) = \int_X \frac{1-e^{-c}}{1 - ye^{-c}}
\prod_{j=1}^{\dim X}\frac{x_j(1 - ye^{-x_j})}{(1-e^{-x_j})}.
\end{eqnarray}

\begin{remark}
The reader is warned of the following possible confusion.
Taking $y=1$ on both sides of (\ref{eqn:chiy}),
one seems to get
\begin{eqnarray*}
\chi(Y) & = & \chi_{-1}(Y) = \int_X \frac{1-e^{-c}}{1 - e^{-c}}
\prod_{j=1}^{\dim X}\frac{x_j(1 - e^{-x_j})}{(1-e^{-x_j})} \\
& = & \int_X \prod_{j=1}^{\dim X} x_j = \chi(X).
\end{eqnarray*}
This is of course absurd.
In fact,
\begin{eqnarray*} \frac{1-e^{-c}}{1 - y e^{-c}}
= \frac{1- e^{-c}}
{(1-y) + y(1 -e^{-c})}
=  \sum_{k \geq 1} \frac{(-y)^{k-1}}{(1-y)^k}(1- e^{-c})^k
\end{eqnarray*}
has a singularity as $y \to 1$,
and one cannot take the limit as above.
\end{remark}

\subsection{Generalization to smooth complete intersection}
\label{sec:CI}

Suppose now $\pi: V \to X$ is a holomorphic vector bundle on $X$ of rank $r$,
$s: X \to V$ a holomorphic section transverse to the zero section.
Then $Y:=s^{-1}(0)$ is a complex submanifold of $X$,
and we have
\begin{eqnarray} \label{eqn:Adjunct2}
&& N_{Y/X} \cong V|_Y.
\end{eqnarray}
Hence from the exact sequence
\begin{eqnarray} \label{eqn:Adjunct3}
0 \to TY \to TX|_Y \to N_{Y/X} \to 0
\end{eqnarray}
by the multiplicative property of $\Ell$,
one gets:
$$\Ell(TY;q,y)\Ell(V|_Y;q;y) = \Ell(TX|_Y;q,y) = \Ell(TX;q,y)|_Y,$$
or equivalently,
$$\Ell(TY;q,y) = \Ell(TX;q,y)\Ell(V;q;y)^{-1}|_Y$$
Hence
$$\chi(Y, \Ell(TY;q,y)) = \chi(Y, \Ell(TX-V;q,y)|_Y)$$
On the other hand,
tensoring the Koszul complex \cite{Gri-Har}
\begin{eqnarray*}
0 & \to & \cO_X(\Lambda^r(V^*)) \xrightarrow{i_s} \cO_X(\Lambda^{r-1}(V^*))
\xrightarrow{i_s}
\dots \\
& \xrightarrow{i_s}&  \cO_X(\Lambda^1(V^*)) \xrightarrow{i_s}  \cO_X
\to \cO_Y \to 0
\end{eqnarray*}
by $\Ell(TX-V;q,y)$,
one gets an exact sequence
\begin{eqnarray*}
0 & \to & \cO_X(\Ell(TX-V;q,y)\Lambda^rV^*)
\xrightarrow{i_s} \cO_X(\Ell(TX-V;q,y)\Lambda^{r-1}(V^*)) \xrightarrow{i_s} \cdots \\
& \xrightarrow{i_s} & \cO_X(\Ell(TX-V;q,y)\Lambda^1V^*)
\xrightarrow{i_s} \cO_X(\Ell(TX-V;q,y)) \\
& \to & \cO_X(\Ell(TX-V;q,y))|_Y \to 0,
\end{eqnarray*}
hence by taking the Riemann-Roch number,
one gets
\begin{eqnarray*}
&&\chi(Y, \Ell(TX-V;q,y)|_Y)
= \chi(X, \Ell(TX-V; q, y)\Lambda_{-1}(V^*)).
\end{eqnarray*}
Therefore,
we have proved the following

\begin{theorem} \label{thm:CI}
Suppose  $Y$ is the zero set of a holomorphic section to a holomorphic
vector bundle $V$ on a compact complex manifold $X$,
which is transverse to the zero section.
Then we have
\begin{eqnarray*}
&& \chi(Y; q, y) = \chi(X, \Ell(TX-V; q, y)\Lambda_{-1}(V^*)).
\end{eqnarray*}
\end{theorem}

Denote by $z_1, \dots, z_r$ the formal Chern roots of $V$.
Then we have by HRR
\begin{eqnarray*}
\chi(Y; q, y)
& = & y^{-\dim Y/2} \int_X \prod_{i=1}^r \prod_{n \geq 1}
\frac{(1-q^{n-1}e^{-z_i})(1-q^ne^{z_i})}{(1-yq^{n-1}e^{-z_i})(1-y^{-1}q^ne^{z_i})} \\
&& \cdot \prod_{j=1}^{\dim X} x_j \prod_{n \geq 1}
\frac{(1-yq^{n-1}e^{-x_j})(1- y^{-1}q^ne^{x_j})}
{(1-q^{n-1}e^{-x_j})(1-q^ne^{x_j})}.
\end{eqnarray*}
Taking $q = 0$,
\begin{eqnarray*}
\chi_{-y}(Y)
& = & \chi\left(X, \Lambda_{-y}(TX)(\Lambda_{-y}V^*)^{-1}\Lambda_{-1}V^*\right) \\
& = & \int_X \prod_{i=1}^r \frac{(1-e^{-z_i})}{(1-ye^{-z_i})}
\cdot \prod_{j=1}^{\dim X} \frac{x_j(1-ye^{-x_j})}{(1-e^{-x_j})}.
\end{eqnarray*}
By (\ref{eqn:Adjunct2}) and (\ref{eqn:Adjunct3}) one can easily deduce
\begin{eqnarray*}
\chi(Y) & = & \int_X \frac{c_r(V)c(TX)}{c(V)} \\
& = & \int_X \prod_{i=1}^r \frac{z_i}{1+z_i}
\cdot \prod_{j=1}^{\dim X}(1+x_j).
\end{eqnarray*}

\subsection{Elliptic genera of projective complete intersections}
\label{sec:General}

Now we are ready to define the elliptic genera of projective complete intersections.
Let $Y\subset \bC\bP^{N-1}$ be a projective variety of codimension $r$,
defined by $r$ homogeneous polynomials $p_1, \dots, p_r$,
of degrees $m_1, \dots, m_r$ respectively.
The virtual bundle
$$(T\bC\bP^{N-1} - \cO_{\bC\bP^{N-1}}(m_1) - \cdots -\cO_{\bC\bP^{N-1}}(m_r))|_Y$$
is the {\em virtual tangent bundle of $Y$ in $\bC\bP^{N-1}$},
and
$$(\cO_{\bC\bP^{N-1}}(m_1) \oplus \cdots  \oplus \cO_{\bC\bP^{N-1}}(m_r))|_Y$$
is the {\em virtual normal bundle of $Y$ in $\bC\bP^{N-1}$}.

The smooth case inspires us to define
$$\chi(Y;q, y) =
\chi(\bC\bP^{N-1}, \Ell(T\bC\bP^{N-1}-V; q,
y)\Lambda_{-1}(V^*)),$$ where $V = \cO(m_1) \oplus \cdots \oplus
\cO(m_k)$. Now Theorem \ref{thm:CI} has the following easy
consequence:

\begin{theorem}
For a nonsingular complete intersection $Y$ in $\bC\bP^{N-1}$, the
elliptic genus of $Y$ defined above coincides with the elliptic
genus of $Y$ as a complex manifold.
\end{theorem}

The relationship between the elliptic genus defined
in this paper and that in \cite{Bor-Lib}
is not clear at present.
Conceivable they should coincide when both are defined.

\section{Elliptic Genera of Complete Intersections in
Projective Spaces as Residues}
\label{sec:Projective}

In this section we express the elliptic genera of complete intersections in projective spaces
as residues.
As motivations we recall the classical calculations for Euler numbers and Hirzebruch
$\chi_y$-genera \cite{Hir}.

\subsection{Euler numbers of smooth hypersurfaces in projective spaces}
Let $Y^N_m$ be a smooth hypersurface of $\bC\bP^{N-1}$ of degree $m$.
Denote the hyperplane class by $H$,
then we have $[Y^N_m] = mH$.
As a warmup,
we recall the computation of the Euler number of $Y^N_m$
(cf. \cite{Hir}).
Since we have an exact sequence
$$0 \to \cO  \to \cO(1)^{\oplus N}
\to T\bC\bP^{N-1} \to 0,
$$
hence we have
$$c(T\bC\bP^{N-1}) = (1+H)^N.$$
Therefore we have
\begin{eqnarray*}
\chi(Y^N_m) & = & \int_{\bC\bP^{N-1}} \frac{mH (1+H)^N}{1 +m H}
= \Res_z \frac{m(1+z)^N}{z^{N-1}(1 +mz)}.
\end{eqnarray*}
Consider the generating series
$\chi_m = \sum_{N = 2}^{\infty} t^{N-1} \chi(Y^N_m)$.
Apply Lemma \ref{lm:Residue} for
$f(z) = z/(1+z)$ one gets first $z = g(t) = t/(1-t)$,
and $f'(z) = 1/(1+z)^2$,
and so
\begin{eqnarray*}
\sum_{N \geq 2} \chi(Y^N_m) t^{N-1}
& = & \frac{mz}{1+m z} (1+z)^2|_{z = \frac{t}{1-t}}
= \frac{mt}{(1-t)^2(1 + (m-1)t)}.
\end{eqnarray*}
From this we deduce the following

\begin{proposition}
Let $Y^N_m$ be a smooth hypersurface of $\bC\bP^{N-1}$ of degree $m$,
then we have
$$\chi(Y^N_m)
= \sum_{k=2}^{N} (-1)^k \begin{pmatrix} N \\ k\end{pmatrix} m^{k-1}.$$
\end{proposition}

\begin{proof}
An elementary calculation gives
\begin{eqnarray*}
&& \frac{mt}{(1-t)^2(1 + (m-1)t)} \\
& = & t\left(\frac{(m-1)^2}{m(1 + (m-1)t)}
+ \frac{m-1}{m(1-t)} + \frac{1}{(1-t)^2} \right)\\
& = & t \left(\frac{(m-1)^2}{m}\sum_{n \geq 0} (1-m)^n t^n
+  \frac{m-1}{m} \sum_{ n \geq 0} t^n
+ \sum_{n \geq 0} (n+1)t^n\right) \\
& = & \sum_{n =1}^{\infty} \left( \frac{(1-m)^{n+1} - (1-(n+1)m)}{m} \right)t^n \\
& = & \sum_{n =1}^{\infty}
\sum_{k=2}^{n+1} (-1)^k \begin{pmatrix} n + 1 \\ k\end{pmatrix} m^{k-1}  t^n.
\end{eqnarray*}
\end{proof}

\subsection{Hirzebruch $\chi_y$-genus of smooth hypersurfaces in projective spaces}
The above calculation can be generalized to the Hirzebruch $\chi_y$-genera \cite{Hir}.
Let $Y^N_m$ be a smooth hypersurface of $\bC\bP^{N-1}$ of degree $m$.
By the multiplicative  properties of the Todd class and $\Lambda_{-y}$,
we have
\begin{eqnarray*}
\td(T\bC\bP^{N-1}) & = & \left(\frac{H}{1-e^{-H}}\right)^{N}, \\
\ch(\Lambda_{-y}(T\bC\bP^{N-1}; q, y)) & = & (1 - ye^{-H})^{N}(1-y)^{-1}.
\end{eqnarray*}
Therefore we have
\begin{eqnarray*}
\chi_{-y}(Y^N_m) & = &  \frac{1}{1-y} \int_{\bC\bP^{N-1}}
\frac{1-e^{-mH}}{1 - ye^{-mH}} (1-ye^{-H})^{N}
\left(\frac{H}{1-e^{-H}}\right)^{N} \\
& = & \frac{1}{1-y} \Res_z \frac{1-e^{-mz}}{1 - ye^{-mz}}
\left(\frac{1-ye^{-z}}{1-e^{-z}}\right)^{N}.
\end{eqnarray*}
We now apply Lemma \ref{lm:Residue} to
$$f(z) = \frac{1-e^{-z}}{1 - ye^{-z}}.$$
First of all,
$$f'(z) = \frac{(1-y)e^{-z}}{(1 - ye^{-z})^2}.$$
Hence $f'(0) = \frac{1}{1-y} \neq 0$.
Secondly, from $f(z) = t$ one finds
$$z = g(t) = - \ln \frac{1-t}{1-yt}.$$
Therefore,
\begin{eqnarray*}
\sum_{N \geq 2} t^{N-1} \chi_{-y}(Y^N_m)
& = &  \frac{1}{1-y} \frac{f(mg(t))}{f'(g(t))} = \frac{1 - e^{-mg(t)}}{1 - y e^{-mg(t)}} \cdot
\frac{(1 - y e^{- g(t)})^2}{(1 - y)e^{-g(t)}} \\
& = &  \frac{1}{1-y} \frac{1 - \left(\frac{1 - t}{1 - y t}\right)^m}
{1 - y \left(\frac{1-t}{1 - yt}\right)^m} \cdot
\frac{\left(1 - y \left(\frac{1-t}{1 - yt}\right)\right)^2}
{(1 - y)\left(\frac{1-t}{1 - yt}\right)} \\
& = &  \frac{(1-yt)^m - (1 - t)^m}{(1-yt)^m - y (1-t)^m}
\cdot \frac{1}{(1-t)(1-yt)}.
\end{eqnarray*}
By Lemma \ref{lm:Residue2} one can also get
$$\chi_{-y}(Y^N_m) = \begin{cases}
\frac{(3m-m^2)(1+y)}{2}, & N = 3, \\
\frac{1+y^2}{6}(m^3-6m^2+12m)+ \frac{y}{2}(2m^3-6m^2+7m), & N = 4.
\end{cases}$$

\subsection{Elliptic genera of hypersurfaces in $\bC\bP^{N-1}$}
\label{sec:Hyp}

Now we generalize to the case of elliptic genera.
Our main result of this section is as follows.

\begin{theorem} \label{thm:EG}
Let $Y^N_m$ be a hypersurface of $\bC\bP^{N-1}$ of degree $m$.
Then we have
\begin{eqnarray}
\chi(Y^N_m; q, y) & = & -\frac{1}{G(y, q)}\oint
\frac{\theta(m z, \tau)}{\theta(v+m z, \tau)}
\left(\frac{\theta(v+z, \tau)}{\theta(z, \tau)}\right)^N dz,
\end{eqnarray}
where
$$G(y,q) = y^{-1/2}\prod_{k \geq 1} \frac{(1-yq^{k-1})(1-y^{-1}q^k)}{(1-q^k)^2}.$$
\end{theorem}

\begin{proof}
By the multiplicative properties of $\Ell$,
we have
\begin{eqnarray*}
\ch(\Ell(T\bC\bP^{N-1}; q, y)) = \frac{1}{G(y, q)}
\left(y^{-1/2}\prod_{k\geq 1} \frac{(1-yq^{k-1}e^{-H})(1- y^{-1}q^ke^{H})}
{(1-q^{k}e^{-H})(1-q^ke^{H})}\right)^{N}.
\end{eqnarray*}
Therefore,
\begin{eqnarray*}
\chi(Y^N_m; q, y) & = & \frac{1}{G(y, q)}
\int_{\bC\bP^{N-1}}y^{1/2}\prod_{k \geq 1}
\frac{(1-q^{k-1}e^{-mH})(1-q^ke^{mH})}{(1-yq^{k-1}e^{-mH})(1-y^{-1}q^ke^{mH})} \\
&& \cdot \left(Hy^{-1/2}\prod_{k\geq 1} \frac{(1-yq^{k-1}e^{-H})(1- y^{-2}q^ke^{H})}
{(1-q^{k-1}e^{-H})(1-q^ke^{H})}\right)^{N} \\
& = & \frac{1}{G(y, q)}\Res_z y^{1/2}\prod_{k \geq 1}
\frac{(1-q^{k-1}e^{-mz})(1-q^ke^{mz})}{(1-yq^{k-1}e^{-mz})(1-y^{-1}q^ke^{mz})} \\
&& \cdot \left(y^{-1/2}\prod_{k \geq 1} \frac{(1-yq^{k-1}e^{-z})(1- y^{-1}q^ke^{z})}
{(1-q^{k-1}e^{-z})(1-q^ke^{z})}\right)^{N}.
\end{eqnarray*}
Now note
$$y^{\frac{1}{2}} \prod_{k \geq 1}
\frac{(1-q^{k-1}e^{-z})(1-q^ke^{z})}{(1-yq^{k-1}e^{-z})(1-y^{-1}q^ke^{z})}
= \frac{\theta(\frac{\sqrt{-1}z}{2\pi}, \tau)}{\theta(v+\frac{\sqrt{-1}z}{2\pi}, \tau)},$$
hence
\begin{eqnarray*}
\chi(Y^N_m; q, y) & = &  \frac{1}{G(y, q)}\Res_z
\frac{\theta(\frac{\sqrt{-1}m z}{2\pi}, \tau)}{\theta(v+\frac{\sqrt{-1}m z}{2\pi}, \tau)}
\left(\frac{\theta(v+\frac{\sqrt{-1}z}{2\pi}, \tau)}
{\theta(\frac{\sqrt{-1}z}{2\pi}, \tau)}
\right)^N \\
& = & \frac{1}{G(y, q)}\frac{1}{2\pi \sqrt{-1}}
\oint \frac{\theta(\frac{\sqrt{-1}m z}{2\pi}, \tau)}{\theta(v+\frac{\sqrt{-1}m z}{2\pi}, \tau)}
\left(\frac{\theta(v+\frac{\sqrt{-1}z}{2\pi}, \tau)}
{\theta(\frac{\sqrt{-1}z}{2\pi}, \tau)}
\right)^N dz \\
& = & -\frac{1}{G(y, q)}\oint
\frac{\theta(m z, \tau)}{\theta(v+m z, \tau)}
\left(\frac{\theta(v+z, \tau)}{\theta(z, \tau)}\right)^N dz.
\end{eqnarray*}
\end{proof}

Introduce a function
$$f(z) =
 \frac{\theta(\frac{\sqrt{-1}z}{2\pi}, \tau)}{\theta(v+\frac{\sqrt{-1}z}{2\pi}, \tau)}.$$
Then we have
\begin{eqnarray*}
\chi(Y^N_m; q, y) & = & \frac{1}{G(y, q)}\Res_z
\frac{f(mz)}{f(z)^{N}}.
\end{eqnarray*}
Note $f(0) = 0$,
and
\begin{eqnarray*}
f'(0)& = & \left. \frac{\sqrt{-1}}{2\pi} \cdot
\frac{\theta'(\frac{\sqrt{-1}z}{2\pi}, \tau)\theta(v + \frac{\sqrt{-1}z}{2\pi}, \tau)
-\theta(\frac{\sqrt{-1}z}{2\pi}, \tau)\theta'(v + \frac{\sqrt{-1}z}{2\pi}, \tau)}
{\theta(v + \frac{\sqrt{-1}z}{2\pi}, \tau)^2}\right|_{z=0} \\
& = & \frac{\sqrt{-1}}{2\pi} \cdot
\frac{\theta'(0, \tau)\theta(v, \tau)-\theta(0, \tau)\theta'(v, \tau)}{\theta(v, \tau)^2}\\
& = & \frac{\sqrt{-1}}{2\pi} \cdot
\frac{\theta'(0, \tau)}{\theta(v, \tau)} \neq 0
\end{eqnarray*}
for $q \neq 0$.
See e.g. \cite{Cha}, p. 65.
Apply Lemma \ref{lm:Residue} for $f$,
we get the following

\begin{corollary}
Suppose $z = g(t)$ is the inverse function of $f(z) =
\frac{\theta(\frac{\sqrt{-1}z}{2\pi}, \tau)}{\theta(v+\frac{\sqrt{-1}z}{2\pi}, \tau)}$
near $z = 0$.
Then we have
\begin{eqnarray*}
&& \sum_{N \geq 2} t^{N-1} \chi(Y^N_m; q, y) = \frac{1}{G(y,
q)}\frac{f(m g(t))}{f'(g(t))}.
\end{eqnarray*}
\end{corollary}

\subsubsection{First few terms}
Since $g(t)$ is hard to find explicitly in this case, we apply
Lemma \ref{lm:Residue2} to get the first few terms. When $N=2$,
$Y^2_m$ consists of $m$ points, hence $\chi(Y^2_m; q, y) = m$. On
the other hand, by Lemma \ref{lm:Residue2} we get
$$\chi(Y^2_m; q, y) = \frac{1}{G(y, q)} \frac{m}{f'(0)},$$
and so we have an identity:
\begin{eqnarray} \label{eqn:f'}
f'(0) = \frac{1}{G(y,q)}.
\end{eqnarray}
Recall
$$G(y,q) = \frac{\theta(v, \tau)}{\sqrt{-1}\eta(q)^3},$$ where $\eta$ is the Dedekind eta
function:
$$\eta(\tau) = q^{\frac{1}{24}}\prod_{k \geq 1} (1 - q^k).$$
Combining with the fact that
$$f'(0) = \frac{\sqrt{-1}}{2\pi}\frac{\theta'(0, \tau)}{\theta(v, \tau)},$$
and (\ref{eqn:f'}),
one recovers the classical result (cf. e.g. \cite{Cha}, p. 71):
$$\theta'(0, \tau) = 2\pi \eta(q)^3.$$
We also get:
\begin{eqnarray} \label{eqn:G}
&& \frac{\theta(v, \tau)}{\theta'(0, \tau)} = \frac{\sqrt{-1}}{2\pi}G(y, q).
\end{eqnarray}

When $N = 3$, we get
$$\chi(Y^3_m; q, y) = \frac{m(m-3)}{2} \frac{f''(0)}{f'(0)}.$$
In particular, when $m = 3$, one gets $\chi(Y^3_3; q, y) = 0$,
this is compatible with the fact that $Y^3_3$ is an elliptic
curve. When $N = 4$,
\begin{eqnarray*}
\chi(Y^4_m; q, y) = \frac{m^3-4m}{6} \frac{f'''(0)}{6f'(0)^3} +
\frac{-2m^2+5m}{2} \frac{f''(0)^2}{f'(0)^5}.
\end{eqnarray*}
Some straightforward calculations yield
\begin{eqnarray*}
f''(0) & = & \left(\frac{\sqrt{-1}}{2\pi}\right)^2
\frac{- 2 \theta'(0, \tau)\theta'(v, \tau)}{\theta(v, \tau)^2}, \\
f'''(0) & = &  \left(\frac{\sqrt{-1}}{2\pi}\right)^3
\left[\frac{\theta'''(0, \tau)}{\theta(v, \tau)}
-3 \frac{\theta'(0, \tau)\theta''(v, \tau)}{\theta(v, \tau)^2}
+ 6\frac{\theta'(0, \tau)\theta'(v, \tau)^2}{\theta(v, \tau)^2}\right].
\end{eqnarray*}
Therefore,
\begin{eqnarray*}
&& \frac{f''(0)}{f'(0)^2} =  \frac{- 2\theta'(v, \tau)}{ \theta'(0, \tau)}
= - \frac{1}{\pi}
\frac{\theta'(v, \tau)}{\eta(q)^3}, \\
&& \frac{f'''(0)}{f'(0)^3} =
\frac{\theta'''(0, \tau)\theta(v, \tau)^2}{\theta'(0, \tau)^3}
-3 \frac{\theta''(v, \tau)\theta(v, \tau)}
{\theta'(0, \tau)^2}
+ 6 \frac{\theta'(v, \tau)^2}{\theta'(0, \tau)^2}.
\end{eqnarray*}

\subsection{Generalization to smooth complete intersections in projective spaces}
Denote by $Y^N_{m_1, \dots, m_r}$ a smooth complete intersection
in $\bC\bP^{N-1}$ defined by $r$ homogeneous polynomials $p_1,
\dots, p_r$ of degrees $m_1, \dots, m_r$ respectively. The above
results can be easily generalized as follows:
\begin{eqnarray*}
&& \sum_{N = r+1}^{\infty} t^{N-1} \chi(Y^N_{m_1, \dots, m_r})
= \frac{1}{(1-t)^2} \prod_{j=1}^r \frac{m_jt}{1 + (m_j - 1)t}, \\
&& \sum_{N = r+1}^{\infty} t^{N-1} \chi_{-y}(Y^N_{m_1, \dots,
m_r}) = \frac{1}{(1-t)(1-yt)}\cdot \prod_{j=1}^r
\frac{(1-yt)^{m_j} - (1 - t)^{m_j}}{(1-yt)^{m_j} - y (1-t)^{m_j}}, \\
&& \sum_{N = r+1}^{\infty} t^{N-1} \chi(Y^N_{m_1, \dots, m_r}; q,
y) = \frac{1}{G(y, q)f'(g(t))}\prod_{j=1}^r f(m_j g(t)),
\end{eqnarray*}
where
$f(z) =
\frac{\theta(\frac{\sqrt{-1}z}{2\pi}, \tau)}{\theta(v+\frac{\sqrt{-1}z}{2\pi}, \tau)}$,
and $g(t)$ is the inverse function of $f$ near $z = 0$.
The first two cases are well-known \cite{Hir}.
In particular,
\begin{eqnarray}
\chi(Y^N_{m_1, \dots, m_r}; q, y) & = & -\frac{1}{G(y, q)}\oint
\prod_{i=1}^r \frac{\theta(m_i z, \tau)}{\theta(v+m_i z, \tau)} \cdot
\left(\frac{\theta(v+z, \tau)}{\theta(z, \tau)}\right)^N dz.
\end{eqnarray}

\section{Generalization to Witten Genera}
\label{sec:Wit}

\subsection{Witten genera}

Let $X$ be a compact spin manifold of real dimension $2k$.
Denote by $T_{\bR}X$ the real tangent bundle and by $\Delta_X$ the spinor bundle of $X$.
Set
\begin{eqnarray*}
&& \Theta_q'(T_{\bR} X)
= \otimes _{n \geq 1} (\Lambda_{q^n}(T_{\bR}X) \otimes S_{q^n}(T_{\bR}X)), \\
&& \Theta_q(T_{\bR} X)
= \otimes _{n \geq 1} (\Lambda_{-q^{n-\frac{1}{2}}}(T_{\bR}X) \otimes S_{q^n}(T_{\bR}X)), \\
&& \Theta_{-q}(T_{\bR} X)
= \otimes _{n \geq 1} (\Lambda_{q^{n-\frac{1}{2}}}(T_{\bR}X) \otimes S_{q^n}(T_{\bR}X)).
\end{eqnarray*}
Witten \cite{Wit} considered the twisted Dirac operators
$d_s \otimes \Theta_q'(T_{\bR}X) = D \otimes \Delta_{T_{\bR}X} \otimes \Theta'_q(T_{\bR}X)$,
$D \otimes  \Theta_q(T_{\bR}X)$ and $D \otimes  \Theta_{-q}(T_{\bR}X)$.
They are related to the four level $1$ representations of loop group $\widetilde{Spin}(2n)$
(cf. \cite{Liu}),
and can be regarded as twisted Dirac operators on the loops space $LX$.
Consider their indices:
\begin{eqnarray*}
&& \sigma_1(LX, q) = q^{-\frac{k}{6}} \Ind d_s \otimes \Theta_q'(T_{\bR}X)
= \int_X \prod_{j=1}^k x_j
\frac{\theta_1(\frac{\sqrt{-1}}{2\pi}x_j, \tau)}{\theta(\frac{\sqrt{-1}}{2\pi}x_j, \tau)}, \\
&& \sigma_2(LX, q) = q^{-\frac{k}{8}} \Ind D \otimes \Theta_q(T_{\bR}X)
= \int_X \prod_{j=1}^k x_j
\frac{\theta_2(\frac{\sqrt{-1}}{2\pi}x_j, \tau)}{\theta(\frac{\sqrt{-1}}{2\pi}x_j, \tau)}, \\
&& \sigma_3(LX, q) = q^{-\frac{k}{8}} \Ind D \otimes \Theta_{-q}(T_{\bR}X)
= \int_X \prod_{j=1}^k x_j
\frac{\theta_3(\frac{\sqrt{-1}}{2\pi}x_j, \tau)}{\theta(\frac{\sqrt{-1}}{2\pi}x_j, \tau)},
 \end{eqnarray*}
where $\{\pm x_j\}$ are the formal Chern roots of $T_{\bR}X \otimes \bC$.
Note the right-hand side of the above indices also make sense when $X$ is not spin.

\subsection{Relationship with elliptic genus}

Now let $X$ be a complex manifold of complex dimension $k$.
One can rewrite (\ref{eqn:chiX}) as
\begin{eqnarray} \label{eqn:chiX2}
&& \chi(X, q, y) = \int_X \prod_{j=1}^d x_j
\frac{\theta(\upsilon + \frac{\sqrt{-1}}{2\pi}x_j, \tau)}{\theta(\frac{\sqrt{-1}}{2\pi}x_j, \tau)}.
\end{eqnarray}
Since by (\ref{eqn:theta1}) and (\ref{eqn:-theta}) we have
\begin{align*}
\theta_{1}(\upsilon,\tau)  & = \theta(\frac{1}{2}+\upsilon,\tau),
\end{align*}
it follows that
\begin{eqnarray} \label{eqn:Sigma}
&& \sigma_1(LX, q) = \chi(X, q, -1),
\end{eqnarray}

\subsection{Relationship with NS elliptic genera}

Assume $X$ is a complex spin manifold.
Then the canonical line bundle $K_X$ has a square root $K_X^{1/2}$.
Define
\begin{eqnarray*}
&& \Ell_{NS}(TX; q, y) \\
& = & (\sqrt{-1}q^{\frac{1}{8}})^{-\dim T} K_X^{1/2} \otimes
\otimes_{n \geq 1}
\left(\Lambda_{-yq^{n-\frac{1}{2}}}(T^*X) \otimes \Lambda_{-y^{-1}q^{n-\frac{1}{2}}}(TX)
\otimes S_{q^n}(T^*X\oplus TX)\right)
\end{eqnarray*}
and
$$\chi_{NS}(X, q, y) = \chi(X, \Ell_{NS}(TX, q, y)).$$
By HRR we have
\begin{eqnarray*}
\chi_{NS}(X, q, y)
& = &  (\sqrt{-1}q^{\frac{1}{8}})^{-k} \int_X \prod_{j=1}^d x_je^{-x_j/2}
\prod_{n \geq 1}\frac{(1-yq^{n-\frac{1}{2}}e^{-x_j})(1- y^{-1}q^{n-\frac{1}{2}}e^{x_j})}
{(1-q^{n-1}e^{-x_j})(1-q^ne^{x_j})} \\
& = & \int_X \prod_j x_j
\frac{\theta_2(\upsilon + \frac{\sqrt{-1}}{2\pi}x_j, \tau)}
{\theta(\frac{\sqrt{-1}}{2\pi}x_j, \tau)}.
\end{eqnarray*}
(Again note the right-hand side makes sense even if $X$ is not spin.)
It follows that
\begin{eqnarray}
&& \sigma_2(LX, q) = \chi_{NS}(X, q, 1), \\
&& \sigma_3(LX, q) = \chi_{NS}(X, q, -1).
\end{eqnarray}

Note when $N$ is even, $\bC\bP^{N-1}$ is spin;
when $m$ is also even,
$Y^N_m$ is spin.
Similar to Theorem \ref{thm:EG},
we have the following:

\begin{theorem} \label{thm:EG2}
Let $Y^N_m$ be a hypersurface of $\bC\bP^{N-1}$ of degree $m$ such that $N-m$ is even.
Then we have
\begin{eqnarray}
\chi_{NS}(Y^N_m; q, y) & = & -\frac{1}{G_{NS}(y, q)}\oint
\frac{\theta(m z, \tau)}{\theta_2(v+m z, \tau)}
\left(\frac{\theta_2(v+z, \tau)}{\theta(z, \tau)}\right)^N dz,
\end{eqnarray}
where
$$G_{NS}(y,q) = (\sqrt{-1}q^{\frac{1}{8}})^{-1}
\prod_{k \geq 1} \frac{(1-yq^{k-\frac{1}{2}})(1-y^{-1}q^{k-\frac{1}{2}})}{(1-q^k)^2}.$$
There is also an easy generalization to the case of complete intersections.
\end{theorem}

It is not hard to see that
\begin{eqnarray} \label{eqn:GNS}
&& \frac{\theta(-\upsilon+\frac{\tau}{2},
\tau)}{\theta_2'(\frac{\tau}{2}, \tau)} = \frac{\sqrt{-1}}{2\pi}
y^{\frac{1}{2}}G_{NS}(y, q).
\end{eqnarray}

\section{Correspondence with  Landau-Ginzburg Orbifolds}
\label{sec:LG}

In this section, we will prove the generalization of the Calabi-Yau/Landau-Ginzburg
correspondence proposed in \cite{Egu-Jin}.
More precisely,
for some hypersurfaces,
the elliptic genera computed geometrically as above
coincide with their algebraic counterparts computed in orbifoldized Landau-Ginzburg models
(see \cite{Egu-Jin, Kaw-Yam-Yan} and the references therein for details).
We will also discuss various generalizations,
including the generalization to the case of NS elliptic genera,
and generalization to complete intersections.

Throughout the rest of the paper,
denote by $C_p$ a small
circle centered at $p$ in the $z$-plane.

\subsection{Elliptic genera case}
\label{sec:EllCase}

\begin{theorem}
\label{thm:ellcase}
Let $Y^{N}_{m}$ be a hypersurface of complex projective space
$\bC\bP^{N-1}$ of degree $m$.
When $(N-m)\upsilon$ is an integer we have
\begin{equation} \label{eqn:EllH}
\chi(Y^{N}_{m};q,y) =
\frac{1}{m}\sum_{a,b=0}^{m-1} e^{2\pi\sqrt{-1}b\upsilon}\left( \frac{\theta
(\frac{m-1}{m}\upsilon+\frac{a}{m}+\frac{b\tau}{m},\tau)}
{\theta(-\frac{1}{m}\upsilon+\frac{a}{m}+\frac{b\tau}{m},\tau)}\right)^{N},
\end{equation}
where $y=e^{2\pi\sqrt{-1}\upsilon}$.
In particular,
when $N=m$, i.e. $Y^{N}_{m}$ is a Calabi-Yau manifold,
(\ref{eqn:EllH}) holds for all $\upsilon$.
\end{theorem}

\begin{proof}
From the distribution of the zeros of theta functions, we can let
$D$ be a fundamental period parallelogram whose boundary,
$C=\partial D$, does not contain the zeros and poles of
$$\frac{\theta(mz,\tau)}{\theta(\upsilon+mz,\tau)}\left(
\frac{\theta(\upsilon+z,\tau)}{\theta(z,\tau)}\right)^{N},$$ and
in this parallelogram it only has simples at
$\frac{a}{m}+\frac{b}{m}\tau-\frac{\upsilon}{m}$, $a,b=0\cdots
m-1$, and another pole at $0$. By Theorem \ref{thm:EG} and residue
theorem, we have
\begin{eqnarray*}
\chi(Y^N_m; q, y) & = & -\frac{1}{G(y, q)}\oint_{C_0}
\frac{\theta(m z, \tau)}{\theta(v+m z, \tau)}
\left(\frac{\theta(v+z, \tau)}{\theta(z, \tau)}\right)^N dz \\
&=&\frac{-1}{G(y,q)}\oint_{C}\frac{\theta(mz,\tau)}
{\theta(\upsilon+mz,\tau)}\left(\frac{\theta(z,\tau)}
{\theta(\upsilon+z,\tau)}\right)^{-N}dz
\\&&+ \frac{1}{G(y,q)}\sum_{a,b=0}^{m-1}\oint_{C_{-\frac{\upsilon}
{m}+\frac{a}{m}+\frac{b}{m}\tau}}\frac{\theta(mz,\tau)}
{\theta(\upsilon+mz,\tau)}\left(\frac{\theta(z,\tau)}
{\theta(\upsilon+z,\tau)}\right)^{-N}dz\\
&=&\frac{1}{G(y,q)}\sum_{a,b=0}^{m-1}
\oint_{C_{-\frac{\upsilon}{m}+\frac{a}{m}+\frac{b}{m}\tau}}\frac{\theta(mz,\tau)}
{\theta(\upsilon+mz,\tau)}\left(\frac{\theta(z,\tau)}
{\theta(\upsilon+z,\tau)}\right)^{-N}dz\\
&=&\frac{1}{mG(y,q)}\sum_{a,b=0}^{m-1} e^{2\pi\sqrt{-1}b\upsilon} G(y,q)\left(
\frac{\theta (\upsilon-\frac{\upsilon}{m}
+\frac{a}{m}+\frac{b}{m}\tau ,\tau)}{\theta(-\frac{\upsilon}{m}+\frac{a}{m}+\frac{b}{m}\tau,\tau)}\right)^{N}\\
&=& \frac{1}{m}\sum_{a,b=0}^{m-1}e^{2\pi\sqrt{-1}b\upsilon}
\left( \frac{\theta(\upsilon-\frac{\upsilon}{m}+\frac{a}{m}+\frac{b}{m}\tau,\tau)}
{\theta(-\frac{\upsilon}{m}+\frac{a}{m}+\frac{b}{m}\tau,\tau)}\right)^{N}.
\end{eqnarray*}
In the above we have used Lemmas \ref{lm:Int1} and \ref{lm:Int2} below.
\end{proof}

We need the following generalization of (\ref{eqn:G}):

\begin{lemma} \label{lm:G}
For $a, b \in \bZ_+$ we have
$$\frac{\theta(-\upsilon+a+b\tau)} {\theta^{'}(a+b\tau)}=
\frac{e^{2\pi\sqrt{-1}b\upsilon}}{2\pi \sqrt{-1}}G(y,q).$$
\end{lemma}

\begin{proof}
By the quasi-periodic properties of theta functions (cf. \S \ref{sec:theta}),
we have
\begin{eqnarray*}
&& \theta(\upsilon+a+b\tau, \tau)
= (-1)^a \theta(\upsilon+b\tau, \tau) \\
& = & (-1)^a (-q^{-\frac{1}{2}}e^{-2\pi\sqrt{-1}(\upsilon+(b-1)\tau)})
\theta(\upsilon+(b-1)\tau, \tau) = \cdots \\
& = & (-1)^{a+b} q^{-\frac{b}{2}} e^{-2\pi\sqrt{-1}b\upsilon}
q^{-\frac{b(b-1)}{2}} \theta(\upsilon, \tau)
= (-1)^{a+b} q^{-\frac{b^2}{2}} e^{-2\pi\sqrt{-1}b\upsilon} \theta(\upsilon, \tau).
\end{eqnarray*}
Hence
\begin{eqnarray*}
&& \theta'(a+b\tau, \tau)
= (-1)^{a+b} q^{-\frac{b^2}{2}} \theta'(0, \tau),
\end{eqnarray*}
and so
\begin{eqnarray*}
&& \frac{\theta(-\upsilon+a+b\tau, \tau)} {\theta^{'}(a+b\tau, \tau)}
= e^{2\pi\sqrt{-1}b\upsilon} \frac{\theta(-\upsilon, \tau)}{\theta'(0, \tau)}
= \frac{e^{2\pi\sqrt{-1}b\upsilon}}{2\pi \sqrt{-1}}G(y,q).
\end{eqnarray*}
In the second equality we have used (\ref{eqn:G}).
\end{proof}

\begin{lemma} \label{lm:Int1}
We have
\begin{eqnarray*}
&&\oint_{C_{\frac{-\upsilon}{m}+\frac{a}{m}+\frac{b}{m}\tau}}\frac{\theta(mz,\tau)}
{\theta(\upsilon+mz,\tau)}\left(
\frac{\theta(\upsilon+z,\tau)}{\theta(z,\tau)}\right)^{N}dz
\\&=&\frac{e^{2\pi\sqrt{-1}b\upsilon}}{m}G(y,q)\left( \frac{\theta
(\upsilon-\frac{\upsilon}{m} +\frac{a}{m}+\frac{b}{m}\tau,\tau)}
{\theta(-\frac{\upsilon}{m}+\frac{a}{m}+\frac{b}{m}\tau,\tau)}\right)^{N}.
\end{eqnarray*}
\end{lemma}

\begin{proof}
We have the following expansions around $z =
-\frac{\upsilon}{m}+\frac{a}{m}+\frac{b}{m}\tau$:
\begin{eqnarray*}
\left(
\frac{\theta(\upsilon+z,\tau)}{\theta(z,\tau)}\right)^{N}
&=& \left( \frac{\theta (\upsilon-\frac{\upsilon}{m}
+\frac{a}{m}+\frac{b}{m}\tau ,\tau)}
{\theta(-\frac{\upsilon}{m}+\frac{a}{m}+\frac{b}{m}\tau,\tau)}\right)^{N}+\cdots,\\
\theta(mz,\tau)&=&\theta(-\upsilon+a+b\tau, \tau)+\cdots,\\
\frac{1}{\theta(\upsilon+mz,\tau)}&=&\frac{1}{m\theta^{'}(a+b\tau,\tau)w}+\cdots,
\end{eqnarray*}
where
$$w = z - (-\frac{\upsilon}{m}+\frac{a}{m}+\frac{b}{m}\tau).$$
Hence
\begin{align*}
&\frac{1}{2\pi
\sqrt{-1}}\oint_{C_{-\frac{\upsilon}{m}+\frac{a}{m}+\frac{b}{m}\tau}}
\frac{\theta(mz,\tau)}{\theta(\upsilon+mz,\tau)}
\left(\frac{\theta(\upsilon+z,\tau)}{\theta(z,\tau)}\right)^{N}dz \\
=&\frac{1}{m}
\frac{\theta(-\upsilon+a+b\tau,\tau)}{\theta^{'}(a+b\tau,\tau)}
\left( \frac{\theta (\upsilon-\frac{\upsilon}{m}
+\frac{a}{m}+\frac{b}{m}\tau,\tau)}
{\theta(-\frac{\upsilon}{m}+\frac{a}{m}+\frac{b}{m}\tau,\tau)}\right)^{N} \\
= & \frac{1}{m}
\frac{e^{2\pi\sqrt{-1}b\upsilon}}{2\pi \sqrt{-1}}G(y,q)
\left( \frac{\theta (\upsilon-\frac{\upsilon}{m}
+\frac{a}{m}+\frac{b}{m}\tau,\tau)}
{\theta(-\frac{\upsilon}{m}+\frac{a}{m}+\frac{b}{m}\tau,\tau)}\right)^{N}.
\end{align*}
In the last equality we have used Lemma \ref{lm:G}.
\end{proof}

\begin{lemma} \label{lm:Int2}
When $(N-m)\upsilon$ is an integer,
$$\oint_{C}\frac{\theta(mz,\tau)}{\theta(\upsilon+mz,\tau)}\left(
\frac{\theta(\upsilon+z\tau,\tau)}{\theta(z,\tau)}\right)^{N}=0.$$
\end{lemma}

\begin{proof}
Use the quasi-periodic properties of $\theta$,
we have
\begin{eqnarray*}
&& \frac{\theta(mz+m,\tau)}{\theta(\upsilon+mz+m,\tau)}\left(
\frac{\theta(\upsilon+z+1,\tau)}{\theta(z+1,\tau)}\right)^{N} \\
& = & \frac{(-1)^m \theta(mz,\tau)}{(-1)^m \theta(\upsilon+mz,\tau)}\left(
\frac{-\theta(\upsilon+z,\tau)}{-\theta(z,\tau)}\right)^{N} \\
& = & \frac{\theta(mz,\tau)}{\theta(\upsilon+mz,\tau)}\left(
\frac{\theta(\upsilon+z,\tau)}{\theta(z,\tau)}\right)^{N}.
\end{eqnarray*}
Furthermore,
\begin{eqnarray*}
&& \frac{\theta(\upsilon+z+\tau,\tau)}{\theta(z+\tau,\tau)}
= e^{-2\pi\sqrt{-1}\upsilon} \frac{\theta(\upsilon+z,\tau)}{\theta(z,\tau)},
\end{eqnarray*}
and so
\begin{eqnarray*}
&& \frac{\theta(\upsilon+mz+m\tau,\tau)}{\theta(mz+m\tau,\tau)} \\
& = & e^{-2\pi\sqrt{-1}\upsilon}
\frac{\theta(\upsilon+mz+(m-1)\tau,\tau)}{\theta(mz+(m-1)\tau,\tau)}
= \cdots \\
& = & e^{-2\pi\sqrt{-1}m \upsilon} \frac{\theta(\upsilon+mz,\tau)}{\theta(mz,\tau)}.
\end{eqnarray*}
Therefore,
when $(N-m)\upsilon$ is an integer, we have
\begin{eqnarray*}
&& \frac{\theta(mz+m\tau,\tau)}{\theta(\upsilon+mz+m\tau,\tau)}\left(
\frac{\theta(\upsilon+z+\tau,\tau)}{\theta(z+\tau,\tau)}\right)^{N}\\
& = & e^{-2\pi\sqrt{-1}(N-m) \upsilon}
\frac{\theta(mz,\tau)}{\theta(\upsilon+mz,\tau)}\left(
\frac{\theta(\upsilon+z,\tau)}{\theta(z,\tau)}\right)^{N} \\
& = & \frac{\theta(mz,\tau)}{\theta(\upsilon+mz,\tau)}\left(
\frac{\theta(\upsilon+z,\tau)}{\theta(z,\tau)}\right)^{N}.
\end{eqnarray*}
It follows that
\begin{eqnarray*}
&&\oint_{C}\frac{\theta(mz,\tau)}{\theta(\upsilon+mz,\tau))}\left(
\frac{\theta(\upsilon+z\tau,\tau)}{\theta(z,\tau)}\right)^{N}dz\\
&=&\int_{p}^{p+1}\frac{\theta(mz,\tau)}{\theta(\upsilon+mz,\tau)}\left(
\frac{\theta(\upsilon+z\tau,\tau)}{\theta(z,\tau)}\right)^{N} dz \\
&& - \int_{p}^{p+1}
\frac{\theta(mz+m\tau,\tau)}{\theta(\upsilon+mz+m\tau,\tau)}\left(
\frac{\theta(\upsilon+z+\tau,\tau)}{\theta(z+\tau,\tau)}\right)^{N}dz\\
&&+\int_{p}^{p+\tau}\frac{\theta(mz,\tau)}{\theta(\upsilon+mz,\tau)}\left(
\frac{\theta(\upsilon+z\tau,\tau)}{\theta(z,\tau)}\right)^{N}dz \\
&& - \int_{p}^{p+\tau}
\frac{\theta(mz+m,\tau)}{\theta(\upsilon+mz+m,\tau)}\left(
\frac{\theta(\upsilon+z+1,\tau)}{\theta(z+1,\tau)}\right)^{N} dz\\
&=&0.
\end{eqnarray*}
\end{proof}

Now we take $\upsilon = \frac{1}{2}$.
When $N-m$ is even,
$(N-m)\upsilon$ is automatically an integer.
So by  (\ref{eqn:Sigma}) and Theorem \ref{thm:ellcase}
we have

\begin{theorem}
\label{thm:signaturecase}
When $N-m$ is an even integer we have
\begin{eqnarray}
&& \sigma_1(LY^N_m, q) = \frac{1}{m}\sum_{a,b=0}^{m-1}(-1)^{b}
\left( \frac{\theta(\frac{1}{2}-\frac{1}{2m}
+\frac{a}{m}+\frac{b}{m}\tau,\tau)}
{\theta(-\frac{1}{2m}+\frac{a}{m}+\frac{b}{m}\tau,\tau)}\right)^{N}.
\end{eqnarray}
\end{theorem}

\subsection{NS elliptic genera version}
One can easily generalize the above results to the case of NS elliptic genera.
First of all we have:

\begin{lemma} \label{lm:GNS}
For $a, b \in \bZ_+$ we have
$$\frac{\theta(-\upsilon+a+b\tau+\frac{\tau}{2})} {\theta'_2(a+b\tau+\frac{\tau}{2})}
= (-1)^a e^{2\pi\sqrt{-1}b\upsilon}
\frac{\theta(-\upsilon+\frac{\tau}{2},
\tau)}{\theta'(\frac{\tau}{2}, \tau)} = (-1)^a
\frac{\sqrt{-1}y^{b+\frac{1}{2}}}{2\pi}G_{NS}(y,q).$$
\end{lemma}

\begin{proof}
By the quasi-periodic properties of theta functions (cf. \S \ref{sec:theta}),
we have
\begin{eqnarray*}
&& \theta(\upsilon+a+b\tau+\frac{\tau}{2}, \tau)
= (-1)^a \theta(\upsilon+b\tau+\frac{\tau}{2}, \tau) \\
& = & (-1)^a (-q^{-\frac{1}{2}}e^{-2\pi\sqrt{-1}(\upsilon+(b-1)\tau+\frac{\tau}{2})})
\theta(\upsilon+(b-1)\tau+\frac{\tau}{2}, \tau) = \cdots \\
& = & (-1)^{a+b} e^{-2\pi\sqrt{-1}b\upsilon}
q^{-\frac{b(b+1)}{2}} \theta(\upsilon+\frac{\tau}{2}, \tau) \\
& = & (-1)^{a+b} e^{-2\pi\sqrt{-1}b\upsilon} q^{-\frac{b(b+1)}{2}}
\theta(\upsilon+\frac{\tau}{2}, \tau).
\end{eqnarray*}
Similarly,
\begin{eqnarray*}
&& \theta_2(\upsilon+a+b\tau+\frac{\tau}{2}, \tau)
= \theta_2(\upsilon+b\tau+\frac{\tau}{2}, \tau) \\
& = & (-q^{-\frac{1}{2}}e^{-2\pi\sqrt{-1}(\upsilon+(b-1)\tau+\frac{\tau}{2})})
\theta_2(\upsilon+(b-1)\tau+\frac{\tau}{2}, \tau) = \cdots \\
& = & (-1)^{b} e^{-2\pi\sqrt{-1}b\upsilon}
q^{-\frac{b(b+1)}{2}} \theta_2(\upsilon+\frac{\tau}{2}, \tau) \\
& = & (-1)^{b} e^{-2\pi\sqrt{-1}b\upsilon} q^{-\frac{b(b+1)}{2}}
\theta_2(\upsilon+\frac{\tau}{2}, \tau).
\end{eqnarray*}
Hence
\begin{eqnarray*}
&& \theta_2'(a+b\tau+\frac{\tau}{2}, \tau)
= (-1)^b q^{-\frac{b(b+1)}{2}} \theta_2'(\frac{\tau}{2}, \tau),
\end{eqnarray*}
and so
\begin{eqnarray*}
\frac{\theta(-\upsilon+a+b\tau+\frac{\tau}{2}, \tau)} {\theta_2^{'}(a+b\tau+\frac{\tau}{2}, \tau)}
& = & (-1)^a e^{2\pi\sqrt{-1}b\upsilon}
\frac{\theta(-\upsilon+\frac{\tau}{2}, \tau)}{\theta_2'(\frac{\tau}{2}, \tau)} \\
& = & (-1)^a \frac{\sqrt{-1}y^{b+\frac{1}{2}}}{2\pi}G_{NS}(y,q).
\end{eqnarray*}
In the second equality we have used (\ref{eqn:GNS}).
\end{proof}

\begin{lemma}
We have
\begin{eqnarray*}
&&\oint_{C_{\frac{-\upsilon}{m}+\frac{a}{m}+\frac{b}{m}\tau+\frac{1}{2m}\tau}}\frac{\theta(mz,\tau)}
{\theta_2(\upsilon+mz,\tau)}\left(
\frac{\theta_2(\upsilon+z,\tau)}{\theta(z,\tau)}\right)^{N}dz
\\&=&\frac{-1}{m}G_{NS}(y,q)(-1)^{a}y^{b+\frac{1}{2}}\left( \frac{\theta_2
(\upsilon-\frac{\upsilon}{m}
+\frac{a}{m}+\frac{b}{m}\tau+\frac{1}{2m}\tau
,\tau)}{\theta(-\frac{\upsilon}{m}+\frac{a}{m}+\frac{b}{m}\tau+\frac{1}{2m}\tau,\tau)}\right)^{N}.
\end{eqnarray*}
\end{lemma}

\begin{proof}
We have the following expansions around $z =
-\frac{\upsilon}{m}+\frac{a}{m}+\frac{b}{m}\tau+ \frac{1}{2m}\tau$:
\begin{eqnarray*}
\left(
\frac{\theta_2(\upsilon+z,\tau)}{\theta(z,\tau)}\right)^{N}
&=& \left( \frac{\theta_2(\upsilon-\frac{\upsilon}{m}
+\frac{a}{m}+\frac{b}{m}\tau+ \frac{1}{2m}\tau,\tau)}
{\theta(-\frac{\upsilon}{m}+\frac{a}{m}+\frac{b}{m}\tau+ \frac{1}{2m}\tau,\tau)}\right)^{N}+\cdots,\\
\theta(mz,\tau)&=&\theta(-\upsilon+a+b\tau+\frac{\tau}{2}, \tau)+\cdots,\\
\frac{1}{\theta_2(\upsilon+mz,\tau)}&=&\frac{1}{m\theta_2^{'}(a+b\tau+ \frac{\tau}{2},\tau)w}+\cdots,
\end{eqnarray*}
where
$$w = z - (-\frac{\upsilon}{m}+\frac{a}{m}+\frac{b}{m}\tau+ \frac{1}{2m}\tau).$$
Hence
\begin{align*}
&\frac{1}{2\pi
\sqrt{-1}}\oint_{C_{-\frac{\upsilon}{m}+\frac{a}{m}+\frac{b}{m}\tau}}
\frac{\theta(mz,\tau)}{\theta_2(\upsilon+mz,\tau)}
\left(\frac{\theta_2(\upsilon+z,\tau)}{\theta(z,\tau)}\right)^{N}dz \\
=&\frac{1}{m}
\frac{\theta(-\upsilon+a+b\tau+ \frac{\tau}{2},\tau)}{\theta_2'(a+b\tau+ \frac{\tau}{2},\tau)}
\left( \frac{\theta_2 (\upsilon-\frac{\upsilon}{m}
+\frac{a}{m}+\frac{b}{m}\tau+ \frac{1}{2m}\tau,\tau)}
{\theta(-\frac{\upsilon}{m}+\frac{a}{m}+\frac{b}{m}\tau+ \frac{1}{2m}\tau,\tau)}\right)^{N} \\
= & \frac{(-1)^a}{m} \frac{
\sqrt{-1}y^{b+\frac{1}{2}}}{2\pi}G_{NS}(y,q) \left( \frac{\theta_2
(\upsilon-\frac{\upsilon}{m} +\frac{a}{m}+\frac{b}{m}\tau+
\frac{1}{2m}\tau,\tau)}
{\theta(-\frac{\upsilon}{m}+\frac{a}{m}+\frac{b}{m}\tau+
\frac{1}{2m}\tau,\tau)}\right)^{N}.
\end{align*}
In the last equality we have used Lemma \ref{lm:GNS}.
\end{proof}

\begin{lemma}
When $N-m$ is an even integer and $(N-m)v$ is
an integer, we have
$$\oint_{C}
\frac{\theta(mz,\tau)} {\theta_2(\upsilon+mz,\tau)}\left(
\frac{\theta_2(\upsilon+z,\tau)}{\theta(z,\tau)}\right)^{N}dz=0,$$
where $C$ is the boundary of a fundamental parallelogram whose
edge does not contain the zeros and poles of
$\frac{\theta(mz,\tau)} {\theta_2(\upsilon+mz,\tau)}\left(
\frac{\theta_2(\upsilon+z,\tau)}{\theta(z,\tau)}\right)^{N}$.
\end{lemma}

\begin{theorem}
\label{thm:Nsgenus} When $N-m$ is even and $(N-m)v$ is an integer,
we have
\begin{eqnarray*}
\chi_{NS}(Y_{m}^{N},q,y)=-\frac{1}{m}
\sum_{a,b=0}^{m-1}(-1)^{a}y^{b+\frac{1}{2}}\left( \frac{\theta_2
(\upsilon-\frac{\upsilon}{m}
+\frac{a}{m}+\frac{b}{m}\tau+\frac{1}{2m}\tau
,\tau)}{\theta(-\frac{\upsilon}{m}+\frac{a}{m}+\frac{b}{m}\tau+\frac{1}{2m}\tau,\tau)}\right)^{N}.
\end{eqnarray*}
\end{theorem}

As a consequence, we have
\begin{corollary}
When $N-m$ is even and $(N-m)v$ is an integer,
we have
\begin{eqnarray*}
\sigma_2(LY_m^N, q) &=&
-\frac{1}{m}\sum_{a,b=0}^{m-1}(-1)^{a}\left( \frac{\theta_2
(\frac{a}{m}+\frac{b}{m}\tau+\frac{1}{2m}\tau ,\tau)}
{\theta(\frac{a}{m}+\frac{b}{m}\tau+\frac{1}{2m}\tau,\tau)}\right)^{N}, \\
\sigma_3(LY_m^N, q) &=&
-\frac{\sqrt{-1}}{m}\sum_{a,b=0}^{m-1}(-1)^{a+b} \left(
\frac{\theta_2(-\frac{1}{2}+\frac{1}{2m}+\frac{a}{m}+\frac{b}{m}\tau+\frac{1}{2m}\tau
,\tau)}
{\theta(\frac{1}{2m}+\frac{a}{m}+\frac{b}{m}\tau+\frac{1}{2m}\tau,\tau)}\right)^{N}
\\&=&-\frac{\sqrt{-1}}{m}\sum_{a,b=0}^{m-1}(-1)^{a+b}\left(
\frac{\theta_3(\frac{1}{2m}+\frac{a}{m}+\frac{b}{m}\tau+\frac{1}{2m}\tau, \tau)}
{\theta(\frac{1}{2m}+\frac{a}{m}+\frac{b}{m}\tau+\frac{1}{2m}\tau,\tau)}\right)^{N}.
\end{eqnarray*}
\end{corollary}

\subsection{Generalization to complete intersections}
Our method admits straightforward generalizations to complete intersections.
Denote $Y^{N}_{m_{1},\cdots,m_{r}}$ is a smooth complete
intersection in $\bC\bP^{N-1}$. Define
$$S(m_{i})=\{\frac{a}{m_{i}}+\frac{b}{m_{i}}\tau|a,b=0, \dots, m_{i}-1\}.$$
In the following we will assume that $\{m_{i}\}$ satisfy
\begin{equation}
S(m_{i})\bigcap S(m_{j})=\phi,\quad i\neq j, \quad i,j=1\cdots r.
\end{equation}
By the same method as above, we can prove
the following Theorems.

\begin{theorem}
When $(N-\sum_{i=1}^{r}m_{i})\upsilon$ is an integer, we have
\begin{eqnarray*}
&&\chi(Y_{m_{1}...m_{r}}^{N};q,y)\\
&=&\sum_{h=1}^{r}\frac{1}{m_{h}}\sum_{a_{h},b_{h}=0}^{m_{h}-1}y^{b_{h}}\left(
\frac{\theta (\upsilon-\frac{\upsilon}{m_{h}}
+\frac{a_{h}}{m_{h}}+\frac{b_{h}}{m_{h}}\tau,\tau
)}{\theta(-\frac{\upsilon}{m_{h}}+\frac{a_{h}}{m_{h}}+\frac{b_{h}}{m_{h}}\tau,\tau)}\right)^{N}
\\&&\cdot\prod_{i\neq h}
\frac{\theta(-\frac{m_{i}\upsilon}{m_{h}}+\frac{m_{i}a_{h}}{m_{h}}
+\frac{m_{i}b_{h}}{m_{h}}\tau,\tau)}
{\theta(\upsilon-\frac{m_{i}\upsilon}{m_{h}}
+\frac{m_{i}a_{h}}{m_{h}}+\frac{m_{i}b_{h}}{m_{h}}\tau,\tau)}.
\end{eqnarray*}
In particular,
\begin{eqnarray*}
&&\sigma_1(LY^{N}_{m_{1}m_{2}...m_{r}}, q) \\
&=&\sum_{h=1}^{r}\frac{1}{m_{h}}\sum_{a_{h},b_{h}=0}^{m_{h}-1}(-1)^{b_h}\left(
\frac{\theta (\frac{1}{2}-\frac{1}{2m_{h}}
+\frac{a_h}{m_{h}}+\frac{b_h}{m_{h}}\tau)}
{\theta(-\frac{1}{2m_{h}}+\frac{a}{m_{h}}+\frac{b}{m_{h}}\tau)}\right)^{N}\\
&&\cdot\prod_{i\neq h}
\frac{\theta(-\frac{m_{i}}{m_{h}}+\frac{m_{i}a_{h}}{m_{h}}+\frac{m_{i}b_{h}}{m_{h}}\tau,\tau)}
{\theta(\frac{1}{2}-\frac{m_{i}}{2m_{h}}+\frac{m_{i}a_{h}}{m_{h}}+\frac{m_{i}b_{h}}{m_{h}}\tau,\tau)}.
\end{eqnarray*}
\end{theorem}

\begin{theorem}
If $(N-\sum_{i=1}^{r}m_{i})$ is an even integer and
$(N-\sum_{i=1}^{r}m_{i})v$ is an integer, then
\begin{eqnarray*}
&&\chi_{NS}(Y_{m_{1}...m_{r}}^{N};q,y)\\
&=&-\sum_{h=1}^{r}\frac{1}{m_{h}}
\sum_{a_h,b_h=0}^{m_h-1}(-1)^{a_h}y^{b_h+\frac{1}{2}}q^{-\frac{1}{2}} \left(
\frac{\theta_2 (\upsilon-\frac{\upsilon}{m_h}
+\frac{a_h}{m_h}+\frac{b_h}{m_h}\tau+\frac{1}{2m_h}\tau ,\tau)}
{\theta(-\frac{\upsilon}{m_h}+\frac{a_h}{m_h}
+\frac{b_h}{m_h}\tau+\frac{1}{2m_h}\tau,\tau)}\right)^{N}
\\&&\cdot\prod_{i\neq h}
\frac{\theta(-\frac{\upsilon m_i}{m_h}+\frac{a_hm_i}{m_h}+\frac{b_hm_i}{m_h}\tau+\frac{m_i}{2m_h}\tau,\tau)}
{\theta_2 (\upsilon-\frac{\upsilon m_i}{m_h}
+\frac{a_hm_i}{m_h}+\frac{b_hm_i}{m_h}\tau+\frac{m_i}{2m_h}\tau, \tau)},
\end{eqnarray*}
and
\begin{eqnarray*}
&&\sigma_2(LY_{m_{1}...m_{r}}^{N}, q)= \chi_{NS}(Y_{m_{1}...m_{r}}^{N}; q, 1),\\
&&\sigma_3(LY_{m_{1}...m_{r}}^{N},q)=\chi_{NS}(LY_{m_{1}...m_{r}}^{N};q,-1).
\end{eqnarray*}
\end{theorem}

\end{document}